# A remark on the Laplacian operator which acts on symmetric tensors

Stepanov S. E., Tsyganok I.I., Aleksandrova I.A.

*Dept. of Math., Finance University under the Government of Russian Federation,*
*125468 Moscow, Leningradsky Prospect, 49-55, Russian Federation*

E-Mail: *s.e.stepanov@mail.ru, i.i.tsyganok@mail.ru, iralex_65@mail.ru*

**Abstract.** More then forty years ago J. H. Samson has defined the Laplacian $\Delta_{sym}$ acting on the space of symmetric covariant *p*-tensors on an *n*-dimensional Riemannian manifold $(M, g)$. This operator is an analogue of the well known Hodge-de Rham Laplacian $\Delta$ which acts on the space of exterior differential *p*-forms $(1 \leq p \leq n)$ on $(M, g)$. In the present paper we will prove that for $n > p = 1$ the operator $\Delta_{sym}$ is the Yano rough Laplacian and show its spectrum properties on a compact Riemannian manifold.

**Key words:** *Riemannian manifold, second order elliptic differential operator on 1-forms, eigenvalues and eigenforms.*

**MSC 2010:** 53C20; 53C21; 53C24

## 1. Definitions and notations

Let $(M, g)$ be a compact oriented $C^\infty$-Riemannian manifold of a dimension $n \geq 2$ with the Levi-Civita connection $\nabla$ and let $S^p M$ be a symmetric tensor product of order $p \geq 1$ of a cotangent bundle $T^*M$ of $M$. On the tensor space $S^p M$ on $M$ we have the *canonical scalar product* $g(\cdot, \cdot)$ and on its $C^\infty$-sections the *global scalar product* $\langle \varphi, \varphi' \rangle = \int_M \frac{1}{p!} g(\varphi, \varphi') dv$ where $dv$ is the volume element of $(M, g)$.

The covariant derivative $\nabla: C^\infty S^p M \to C^\infty(T^*M \otimes S^p M)$ has the formal adjoint operator $\delta = \nabla^*: C^\infty(T^*M \otimes S^p M) \to C^\infty S^p M$ which is uniquely defined by the formula $\langle \nabla \cdot, \cdot \rangle = \langle \cdot, \nabla^* \cdot \rangle$ (see [1, p. 460]). Furthermore we can define (see also [1, p. 514]) the operator $\delta^*: C^\infty S^p M \to C^\infty S^{p+1} M$ which is the adjoint operator of $\delta: C^\infty S^{p+1} M \to C^\infty S^p M$ with respect to the global product $\langle \cdot, \cdot \rangle$.

More then forty years ago J. H. Samson has defined (see [2]) the Laplacian operator $\Delta_{sym} = \delta \delta^* - \delta^* \delta: C^\infty S^p M \to C^\infty S^p M$. This operator is an analogue of the well known Hodge-de Rham Laplacian $\Delta: C^\infty \Lambda^p M \to C^\infty \Lambda^p M$ which acts on $C^\infty$-sections of the bundle $\Lambda^p M$ of covariant skew-symmetric tensors of degree $p$ $(1 \leq p \leq n)$ on $M$ and is defined by $\Delta = d\delta + \delta d$ for the exterior differential $d: \Lambda^p M \to \Lambda^{p+1} M$ (see [1, p. 34]).

The operator $\Delta_{sym}$ is studied in the following papers [2]; [3]; [4]; [5] and [6].

This paper is organized as follows. The next section summarizes the basic properties of $\Delta_{sym}$: $C^\infty S^p M \to C^\infty S^p M$ for the case $p = 1$. Section with the number three expresses our results on infinitesimal conformal and projective transformations. The fourth section of the present paper shows spectrum properties of $\Delta_{sym}$ on an *n*-dimensional compact Riemannian manifold for the case $n > p = 1$. And in the last section we prove a theorem about eigenvalues of the Hodge-de Rham Laplacian $\Delta$ which acts on closed 1-forms.

## 2. The Yano rough Laplacian

We proved in [6] that for $p = 1$ the Weitzenböck decomposition formula for $\Delta_{sym} = \delta\delta^* - \delta^*\delta$ has the form $\Delta_{sym} = \delta\nabla - Ric$ where *Ric* is the Ricci tensor of (*M, g*) and $\delta\nabla$ the *Bochner rough Laplacian* which is also denoted by $\nabla^*\nabla$ (see [1, p. 54]). Next, thanks to the well-known Weitzenböck decomposition formula $\Delta = \delta\nabla + Ric$ for the Hodge-de Rham Laplacian $\Delta: C^\infty T^*M \to C^\infty T^*M$ we concluded that $\Delta_{sym} = \Delta - 2\,Ric$. After that, using the equation $\Delta_{sym} = \Delta - 2\,Ric$ we can define the differential operator $\square: C^\infty TM \to C^\infty TM$ such that $\square = \Delta - 2\,Ric^*$ for the linear symmetric operator $Ric^*$ which is associated with the Ricci tensor *Ric* and defined by the identity $Ric(X,Y) = g(Ric^* X, Y)$ for any $X, Y \in C^\infty TM$ (see also [8, p. 40]). In turn, we recall that more then forty years ago the operator $\square$ was used by K. Yano (see [8]) for the investigation of local isometric, conformal, affine and projective transformations of compact Riemannian manifolds. Based on the above, we will call $\Delta_{sym}$ the *Yano rough Laplacian* when $p = 1$. Hence, the following proposition is true.

**Lemma.** *Let* (*M, g*) *be an n-dimensional* $(n \geq 2)$ *Riemannian manifold. For $p = 1$ the Samson Laplacian* $\Delta_{sym}: C^\infty S^p M \to C^\infty S^p M$ *is the Yano rough Laplacian.*

We recall here that the vector field $\xi$ on (*M, g*) is called an *infinitesimal harmonic transformation* if the one-parameter group of infinitesimal point transformations of (*M, g*) generated by $\xi$ consists of harmonic diffeomorphisms (see [6]). In turn, we have proved in [6] that the vector field $\xi$ is an infinitesimal harmonic transformation on (*M, g*) if and only if $\Delta_{sym}\omega = 0$ for the 1-form $\omega$ dual to the vector field $\xi$ with respect to the metric *g*, i.e. $\omega(X) = g(\xi, X)$ for an arbitrary vector field $X \in C^\infty TM$. In this case, we adopt the following notation $\xi := \omega^\#$.

In particular, *holomorphic vector fields* on nearly Kählerian manifolds (see [9]) and vector fields that transform a Riemannian metrics into Ricci soliton metrics (see [9]) are examples of infinitesimal harmonic transformations. Therefore, all forms which are dual to these vector fields belong to $Ker\,\Delta_{sym}$.

On the other hand, a vector field $\xi$ is called a *Killing vector field* or, in other words an *infinitesimal isometric transformation* if the one-parameter group of infinitesimal transformations of $(M, g)$ generated by $\xi$ consists of isometric diffeomorphisms. An arbitrary Killing vector field $\xi$ satisfies the condition $\delta^*\omega = 0$ where $\xi := \omega^{\#}$. On the other hand, according to the Yano's theorem (see [8, p. 44]; [10]) a vector field $\xi$ on a compact Riemannian manifold $(M, g)$ is a Killing vector field if and only if $\Delta_{sym}\omega = 0$ and $\delta\omega = 0$. The vector space of 1-forms dual to globally defined Killing vector fields has the finite dimension $k_1(M) \leq \frac{1}{2} n (n + 1)$. The dimension $k_1(M)$ has been named the *first Killing number*. Moreover, we have proved in [7] that the number $k_1(M)$ is a scalar projective invariant of $(M, g)$.

### 3. Conformal Killing and projective Killing 1-forms

A real number $\lambda$, for which there is a form $\omega \in C^{\infty}T^*M$ (not identically zero) such that $\Delta_{sym}\omega = \lambda\omega$, is called an *eigenvalue* of $\Delta_{sym}$ and the corresponding $\omega \in C^{\infty}T^*M$ is called an *eigenform* of $\Delta_{sym}$ corresponding to $\lambda$. Next, we consider two examples of eigenforms of $\Delta_{sym}$.

*Conformal Killing vector fields* can be considered as a natural generalization of Killing vector fields. They are also called *infinitesimal conformal transformations* because any conformal Killing vector $\xi$ generates a local one-parameter group of conformal diffeomorphisms of $(M, g)$.

Consider an *n*-dimensional compact orientable Riemannian manifold $(M, g)$. Lichnerowicz has shown (see [8, p. 47]) that a necessary and sufficient condition for $\xi$ to be a *conformal Killing vector field* on $(M, g)$ is

$$\Delta_{sym}\omega + (1 - 2/n)\delta^*\delta\omega = 0 \qquad (3.1)$$

for the 1-form $\omega$ dual to the vector field $\xi$ with respect to the metric $g$. This 1-form is called *conformal Killing form* (see, for example, [7]). Let the eigenform $\omega$ of $\Delta_{sym}$ be a conformal Killing form on an *n*-dimensional ($n > 2$) compact and oriented Riemannian manifold $(M, g)$ then

$$\lambda\langle\omega,\omega\rangle = -n^{-1}(n-2)\langle\omega,\delta^*\delta\omega\rangle = -n^{-1}(n-2)\langle\delta\omega,\delta\omega\rangle.$$

From these equations, we deduce the following inequality

$$\lambda = -(1 - 2/n)\frac{\langle\delta\omega,\delta\omega\rangle}{\langle\omega,\omega\rangle} \leq 0.$$

For the second example we consider a *projective Killing vector field* or, in other words an *infinitesimal projective transformation* (see [8, p. 45]) which satisfies the equation $\Delta_{sym}\omega = 2(n+1)^{-1}\delta^*\delta\omega$ for the form $\omega$ dual to $\xi$. This 1-form will be called *projective Killing form*. Let the eigenform $\omega$ of $\Delta_{sym}$ be a projective Killing form on a compact and oriented Riemannian manifold $(M, g)$. In this case, we have

$$\lambda \langle \omega, \omega \rangle = 2(n+1)^{-1} \langle \omega, \delta^* \delta \omega \rangle = 2(n+1)^{-1} \langle \delta \omega, \delta \omega \rangle$$

and consequently the following inequality holds

$$\lambda = 2(n+1)^{-1} \frac{\langle \delta \omega, \delta \omega \rangle}{\langle \omega, \omega \rangle} \geq 0.$$

## 4. Spectral properties of the Yano rough Laplacian

We recall that all nonzero eigenforms corresponding to a fixed eigenvalue $\lambda$ form a vector subspace of $C^\infty T^* M$ denoted by $V_\lambda(M)$ and called the eigenspace corresponding to the eigenvalue $\lambda$.

The following theorem about eigenvalues of $\Delta_{sym}$ and their corresponding forms is valid.

**Theorem 2.** *Let (M, g) be an n-dimensional $(n \geq 2)$ compact and oriented Riemannian manifold and $\Delta_{sym}: C^\infty T^* M \to C^\infty T^* M$ be the Yano rough Laplacian.*

1) *Suppose the Ricci tensor is negative then an arbitrary eigenvalue $\lambda$ of $\Delta_{sym}$ is positive.*

2) *The eigenspaces of $\Delta_{sym}$ are finite dimensional.*

3) *The eigenforms corresponding to distinct eigenvalues are orthogonal.*

**Proof.** 1) Let $\varphi \in V_\lambda(M)$ be a non-zero eigentensor corresponding to the eigenvalue $\lambda$, that is $\Delta_{sym} \omega = \lambda \omega$ then we can rewrite the formula $\Delta_{sym} = \delta \nabla - Ric$ in the form

$$\lambda \langle \omega, \omega \rangle = -\int_M Ric(\xi, \xi) dv + \langle \nabla \omega, \nabla \omega \rangle. \tag{4.1}$$

where $\xi$ is the vector field dual to the 1-form $\omega$. If we suppose that the Ricci tensor is negative and we denote by $-r$ the largest (negative) eigenvalue of matrix $\|Ric\|$ on (M, g) then $Ric(\xi, \xi) \leq -r g(\xi, \xi)$. In this case from the inequality (4.1), we conclude

$$\lambda \langle \omega, \omega \rangle \geq r \langle \omega, \omega \rangle + \langle \nabla \omega, \nabla \omega \rangle > 0.$$

2) The eigenspaces of $\Delta_{sym}$ are finite dimensional because $\Delta_{sym}$ is an elliptic operator.

3) Let $\lambda_1 \neq \lambda_2$ and $\omega_1, \omega_2$ be the corresponding eigenforms. Then $\langle \Delta_{sym} \omega_1, \omega_2 \rangle = \lambda_1 \langle \omega_1, \omega_2 \rangle$ and $\langle \Delta_{sym} \omega_1, \omega_2 \rangle = \langle \omega_1, \lambda_2 \omega_2 \rangle = \lambda_2 \langle \omega_1, \omega_2 \rangle$. Therefore $0 = (\lambda_1 - \lambda_2) \langle \omega_1, \omega_2 \rangle$ and since $\lambda_1 \neq \lambda_2$ it follows that $\langle \omega_1, \omega_2 \rangle = 0$, that is, $\omega_1$ and $\omega_2$ are orthogonal.

In particular, for the case $n = 2$ we have the following theorem.

**Theorem 3**. *Let (M, g) be a 2-dimensional compact and oriented Riemannian manifold. Then the first eigenvalue $\lambda_1$ of the Yano rough Laplacian $\Delta_{sym}: C^\infty T^* M \to C^\infty T^* M$ is a non-negative number.*

**Proof.** We compute that $g(\delta^*\omega, \delta^*\omega) \geq 4\,n^{-1}(\delta\omega)^2$ for any $\omega \in C^\infty T^*M$. This elementary algebraic fact can be rewritten as $2^{-1} g(\delta^*\omega, \delta^*\omega) - (\delta\omega)^2 \geq -n^{-1}(n-2)(\delta\omega)^2$. Integration by parts yields the following integral inequality $\langle \Delta_{sym}\,\omega, \omega \rangle \geq -n^{-1}(n-2) \int_M (\delta\omega)^2 dv$ where the operator $\Delta_{sym}$ satisfies the identity $\langle \Delta_{sym}\,\omega, \omega \rangle = \langle \delta^*\omega, \delta^*\omega \rangle - \langle \delta\,\omega, \delta\,\omega \rangle$, which follows immediately from its definition. The inequality proves our theorem.

We consider now the $n$-dimensional ($n \geq 2$) Einstein manifold $(M, g)$ where $Ric = \dfrac{s}{n} g$ and $s$ is a constant (see [1, p. 44]). In this case we can rewrite the formula $\Delta_{sym} = \Delta - 2Ric$ in the form

$$\Delta_{sym} = \Delta - 2\frac{s}{n} g. \tag{4.2}$$

From (4.2) we conclude that the following theorem is true.

**Theorem 4.** *Let $(M, g)$ be an n-dimensional ($n \geq 2$) compact and oriented Einstein manifold $(M, g)$ then*
1) *if $s > 0$ then any 1-form which is dual to an infinitesimal harmonic transformation is an eigenform of $\Delta$ corresponding to the eigenvalue $2\dfrac{s}{n}$ and the converse is also true;*
2) *if $s < 0$ then any harmonic 1-form is an eigenform of $\Delta_{sym}$ corresponding to the eigenvalue $-2\dfrac{s}{n}$ and the converse is also true.*

Using the general theory of elliptic operators on a compact $(M, g)$ it can be proved that $\Delta_{sym}$ has a discrete spectrum, denoted by $Spec\,\Delta_{sym}$, consisting of real eigenvalues of finite multiplicity which accumulate only at infinity. In symbols, we have $Spec\,\Delta_{sym} = \{0 \leq |\lambda_1| \leq |\lambda_2| \leq \ldots \to +\infty\}$. In addition, if we suppose that the Ricci tensor $Ric$ is negative then $Spec\,\Delta_{sym} = \{0 < \lambda_1 \leq \lambda_2 \leq \ldots \to +\infty\}$. Moreover, here we have the following:

**Theorem 5.** *Let $(M, g)$ be an n-dimensional ($n \geq 2$) compact and oriented Riemannian manifold. Suppose the Ricci tensor Ric is negative, then the first eigenvalue $\lambda_1$ of the Yano rough Laplacian $\Delta_{sym}$: $C^\infty T^*M \to C^\infty T^*M$ satisfies the inequality $\lambda_1 \geq 2r$ for the largest (negative) eigenvalue $-r$ of matrix $\|Ric\|$ on $(M, g)$. The equality $\lambda_1 = 2r$ is attained for some harmonic eigenform $\omega \in C^\infty T^*M$ and in this case the multiplicity of $\lambda_1$ is less than or equals to the Betti number $b_1(M)$.*

**Proof.** Let $(M, g)$ be an $n$-dimensional compact and oriented Riemannian manifold. Suppose that the Ricci tensor is negative. Denote by $-r$ the largest (negative) eigenvalue of matrix $\|Ric\|$. Then from the formula $\Delta_{sym} = \Delta - 2Ric$ we obtain the inequality

$$\langle \Delta_{sym}\,\omega, \omega \rangle \geq 2r\langle \omega, \omega \rangle + \langle \Delta\omega, \omega \rangle \tag{4.3}$$

for any $\omega \in T^*M$. Then for an eigenform $\omega$ corresponding to an eigenvalue $\lambda$, (4.4) becomes the inequalities

$$\lambda \langle \omega, \omega \rangle \geq 2r\langle \omega, \omega \rangle + \langle \Delta\omega, \omega \rangle \geq 2r\langle \omega, \omega \rangle \tag{4.4}$$

which prove that

$$\lambda_1 \geq 2r > 0. \tag{4.5}$$

If the equality is valid in (4.5), then from (4.4) we obtain $\Delta\omega = 0$. In this case $\omega$ is a harmonic 1-form and, so the multiplicity of $\lambda_1$ is less than or equals to the Betti number $b_1(M)$ because the number of linearly independent (with constant real coefficients) harmonic 1-forms on $(M, g)$ is equal to the Betti number $b_1(M)$ of $(M, g)$ (see [11]). The proof is complete.

Suppose now that $(\mathbb{H}^n, g_0)$ is a compact $n$-dimensional hyperbolic manifold with standard metric $g_0$ having constant sectional curvature equal to $-1$. In this case, from the theorem above we obtain the following corollary.

**Corollary**. Let $(\mathbb{H}^n, g_0)$ *be an n-dimensional compact and oriented hyperbolic manifold then the first eigenvalue* $\lambda_1$ *of the Yano rough Laplacian* $\Delta_{sym}: C^\infty T^*M \to C^\infty T^*M$ *satisfies the inequality* $\lambda_1 \geq 2$. *The equality* $\lambda_1 = 2$ *is attained if and only if* $n = 2$. *In this case the multiplicity of* $\lambda_1$ *is equal to the Betti number* $b_1(\mathbb{H}^2)$.

**Proof**. Let $(M, g)$ be a compact and oriented model of hyperbolic space $(\mathbb{H}^n, g_0)$ with standard metric $g_0$ having constant sectional curvature equal to $-1$ then $\lambda_1 \geq 2$. At the same time it is well known (see [12]) that $L^2$-harmonic $p$-forms appear on a simply connected complete hyperbolic manifold $(M, g)$ of constant sectional curvature $-1$ if and only if $n = 2p$. Therefore, if $(M, g)$ is a compact and oriented model of hyperbolic space $(\mathbb{H}^n, g_0)$ then the equality $\lambda_1 = 2$ is attained if and only if $n = 2$. In this case the multiplicity of $\lambda^r$ is equal to the Betti number $b_1(\mathbb{H}^2)$.

## 5. Appendix

Finally, we prove the following theorem which is dual to the above theorem with the number 5.

**Theorem 6**. *Let $(M, g)$ be an $n$-dimensional ($n \geq 2$) compact and oriented Riemannian manifold and $\mu_1$ be a first eigenvalue of the Laplacian $\Delta: C^\infty T^*M \to C^\infty T^*M$ such that the corresponding 1-form $\omega \in C^\infty T^*M$ is a coclosed form. Moreover, suppose that the Ricci tensor Ric is positive, then $\mu_1 \geq 2\rho$ for the smallest (positive) eigenvalue $\rho$ of matrix $\|Ric\|$ on $(M, g)$. The equality $\mu_1 = 2\rho$ is attained for* some Killing eigenform $\omega \in C^\infty T^*M$ and the multiplicity of $\mu_1$ is less than or equals to the Killing number $k_1(M)$.

**Proof**. Let $\omega$ be a coclosed eigenvalue form of $\Delta$ corresponding to an eigenvalue $\mu$ of $\Delta$ then from the formula $\Delta_{sym} = \Delta - 2Ric$ we obtain the integral equality

$$\mu \langle \omega, \omega \rangle = \langle \delta^*\omega, \delta^*\omega \rangle + 2\int_M Ric(\xi,\xi) dv \qquad (5.1)$$

where $\xi := \omega^{\#}$. Now, if we assume that the Ricci tensor is positive and denote by $\rho$ the smallest (positive) eigenvalue of the matrix $\|Ric\|$, then we have $Ric(X,X) \geq \rho g(X,X)$ for an arbitrary vector field $X \in C^\infty TM$. In this case, thanks to (5.1), we have $\mu_1 \geq 2\rho$. On the other hand, if $\mu_1 = 2\rho$ then from (4.6) we conclude that $\delta^*\omega = 0$. Hence $\xi := \omega^{\#}$ is a Killing vector field. The theorem is proved.